\documentclass[a4paper,10pt]{article}
\usepackage{array}

\usepackage{graphicx}
\usepackage{amsmath}
\usepackage{bbm}
\usepackage[sort&compress,numbers]{natbib}
\usepackage{amsfonts}
\usepackage{amssymb}
\usepackage{amsmath}
\usepackage{latexsym}
\usepackage{color}
\usepackage{url}
\usepackage{ntheorem}
\usepackage{braket}
\usepackage[normalem]{ulem}
\usepackage{tikz}
\usepackage{enumerate}
\usepackage[colorlinks=true]{hyperref}
\newcommand{\ringg}{{\red{\mathring{g}}}}% this is the background metric, equal to the hyperbolic one near the conformal boundary at infinity, but not necessarily everywhere
\newcommand{\hypg}{{\red{b}}}%{{\ringg}}% this is the hyperbolic metric

\newcommand{\enmom}{\red{m}}
\newcommand{\benmom}{{\mathbf{m}}}

\newcommand{\myone}{\red{1}}
\newcommand{\mytwo}{\red{2}}
\newcommand{\myi}{\red{i}}

\newcommand{\Hyp}{{\mathbb H}}

\newcommand\ben{\begin{enumerate}}
\newcommand\een{\end{enumerate}}
\newcommand\bit{\begin{itemize}}
\newcommand\eit{\end{itemize}}

\newcommand{\qedskip}{\qed\medskip}

\newcommand{\ringh}{\,\mathring{\! h}{}}

\newcommand{\red}[1]{{\color{red}#1}}

\newcounter{mnotecount}[section]

\renewcommand{\themnotecount}{\thesection.\arabic{mnotecount}}

\newcommand{\mnote}[1]%{}
{\protect{\stepcounter{mnotecount}}$^{\mbox{\footnotesize
$%\!\!\!\!\!\!\,
\bullet$\themnotecount}}$ \marginpar{%\color{red}%
\raggedright\tiny\em
$\!\!\!\!\!\!\,\bullet$\themnotecount: #1} }

\newcommand{\ptcr}[1]{{\color{red}\mnote{{\color{red}{\bf ptc:}
#1} }}}

\newtheorem{theorem}{\sc  Theorem\rm}[section]

\newtheorem{thm}[theorem]{\sc  Theorem\rm}

\newtheorem{conjecture}[theorem]{\sc  Conjecture\rm}

\newtheorem{remark}[theorem]{\sc Remark\rm}
\newtheorem{Remark}[theorem]{\sc Remark\rm}

\newcommand{\jlcax}[1]{}
\newcommand{\eean}{\nonumber\end{eqnarray}}

\newcommand{\kk}[1]{}%{\mnote{{\bf If we consider the KK case:} #1}}

\newcommand{\mcH}{{\mycal H}}

\newcommand{\beq}{\begin{equation}}

\newcommand{\FS}       %{F_1} %
                  {F}
\newcommand{\HS} %{F_2}
       {H_{\mbox{\scriptsize volume}}}

{\ptc{this should be removed in the oberwolfach version}}%

\newcommand{\zD}{\mathring{D}}%
\newcommand{\eeal}[1]{\label{#1}\end{eqnarray}}
\newcommand{\bed}{\begin{deqarr}}
\newcommand{\eed}{\end{deqarr}}
\newcommand{\bedl}[1]{\begin{deqarr}\label{#1}}
\newcommand{\eedl}[2]{\arrlabel{#1}\label{#2}\end{deqarr}}

\newcommand{\mcU}{{\mycal U}}

\newcommand{\bel}[1]{\begin{equation}\label{#1}}
\newcommand{\bea}{\begin{eqnarray}}
\newcommand{\bean}{\begin{eqnarray}\nonumber}
\newcommand{\beal}[1]{\begin{eqnarray}\label{#1}}
\newcommand{\eea}{\end{eqnarray}}

\newcommand{\qed}{\hfill $\Box$ \medskip}
\newcommand{\proof}{\noindent {\sc Proof:\ }}
\newcommand{\be}{\begin{equation}}
\newcommand{\eeq}{\end{equation}}
\newcommand{\ee}{\end{equation}}
\newcommand{\beqa}{\begin{eqnarray}}
\newcommand{\eeqa}{\end{eqnarray}}
\newcommand{\beqan}{\begin{eqnarray*}}
\newcommand{\eeqan}{\end{eqnarray*}}
\newcommand{\ba}{\begin{array}}
\newcommand{\ea}{\end{array}}

\newcommand{\Id}{\mbox{\rm Id}} %identity matrix
\newcommand{\mcV}{{\mycal V}}

\newcommand{\warn}[1]%{}%{}
{\protect{\stepcounter{mnotecount}}$^{\mbox{\footnotesize
$%\!\!\!\!\!\!\,
\bullet$\themnotecount}}$ \marginpar{%\color{red}%
\raggedright\tiny\em
$\!\!\!\!\!\!\,\bullet$\themnotecount: {\bf Warning:} #1} }

\newcommand{\R}{\mathbb R}

\newcommand{\N}{\mathbb N}

\newcommand{\eq}[1]{(\ref{#1})}

\newcommand{\ptc}[1]{\mnote{{\bf ptc:}#1}}

\newcommand{\beqar}{\begin{deqarr}}
\newcommand{\eeqar}{\end{deqarr}}

\newcommand{\beaa}{\begin{eqnarray*}}
\newcommand{\eeaa}{\end{eqnarray*}}

\newcommand{\tr}{\mbox{tr}}
\newcommand{\zg}{\mathring{g}}

\newcommand{\znabla}{\mathring{\nabla}}

\definecolor{applegreen}{rgb}{0.55, 0.71, 0.0}
\definecolor{armygreen}{rgb}{0.29, 0.33, 0.13}

\definecolor{caribbeangreen}{rgb}{0.0, 0.8, 0.6}

\definecolor{orange}{RGB}{255,127,0}

\newcommand{\purple}[1]{{\color{purple} #1}}

\newcommand{\notjdg}[1]{}
\renewcommand{\notjdg}[1]{\protect\ptcr{purple text only in the arxiv version, not for jdg}\purple{#1}} % this is for the arxiv version
\newcommand{\onlyjdg}[1]{\mnote{green only in jdg}\green{#1}}
\newcommand{\nojdg}[1]{\ptcr{purple only in arxiv, not in jdg}\purple{#1}}

\renewcommand{\myone}{1}
\renewcommand{\mytwo}{2}
\renewcommand{\myi}{i}

\newcommand{\normp}{|}

\DeclareFontFamily{OT1}{rsfs}{}
\DeclareFontShape{OT1}{rsfs}{CGNPm}{n}{ <-7> rsfs5 <7-10> rsfs7 <10-> rsfs10}{}
\DeclareMathAlphabet{\mycal}{OT1}{rsfs}{CGNPm}{n}

{\catcode `\@=11 \global\let\AddToReset=\@addtoreset}
\AddToReset{equation}{section}

{\catcode `\@=11 \global\let\AddToReset=\@addtoreset}
\AddToReset{figure}{section}

{\catcode `\@=11 \global\let\AddToReset=\@addtoreset}
\AddToReset{table}{section}

\renewcommand{\red}[1]{#1}
\renewcommand{\nojdg}[1]{}
\renewcommand{\notjdg}[1]{}
\renewcommand{\onlyjdg}[1]{#1}
\begin{document}
\title{The hyperbolic positive energy theorem\protect\thanks{Preprint UWThPh-2019-3} }

\author{Piotr T. Chru\'{s}ciel\thanks{Beijing Institute of Mathematical Sciences and Applications, Huairou, and Center for Theoretical Physics of the Polish Academy of Sciences, Warsaw} \thanks{
{\sc Email} \protect\url{piotr.chrusciel@univie.ac.at}, {\sc URL} \protect\url{homepage.univie.ac.at/piotr.chrusciel}}
\\
%\author[Delay]
{Erwann
Delay}\thanks{Universit\'e d'Avignon, Laboratoire de Math\'ematiques d'Avignon (EA 2151), 301 rue Baruch de Spinoza,
F-84916 Avignon, France} \thanks{
{\sc Email} \protect\url{erwann.delay@univ-avignon.fr}, {\sc URL} \protect\url{https://erwanndelay.wordpress.com/}}
}
\maketitle

\begin{abstract}
We show that the causal-future-directed character of the energy-momentum vector of $n$-dimensional asymptotically hyperbolic  Riemannian manifolds with spherical conformal infinity, $n\ge 3$, can be traced back to that of  asymptotically Euclidean general-relativistic initial data sets satisfying the dominant energy condition.
\end{abstract}

\tableofcontents
%----------------------------------------------------------------------------%
\section{Introduction}

An interesting global invariant of asymptotically hyperbolic general relativistic initial data sets is the energy-momentum vector $\benmom \equiv (\enmom_\mu)$
 \cite{ChHerzlich,Wang,ChNagyATMP}
 (compare~\cite{AbbottDeser,ChruscielSimon}).
 In the case of spherical conformal infinity, it is known that $\benmom$ is timelike future pointing under a spin condition~\cite{ChHerzlich,Wang,CJL,Maerten,GHHP}.
 The object of this work is to show how to remove this condition.

Our analysis will rely on the following  rigidity conjecture:

\begin{conjecture}
\label{T12XII18.1}
Let  $(M,g,K )$ be an $n$-dimensional general relativistic initial data set, $n\ge 3$, satisfying the dominant energy condition such that
$g$ is flat outside of a compact set and $K$ vanishes there.
 Then $(M,g)$ can be isometrically embedded in Minkowski space-time so that $K$ is the second fundamental form of the embedded hypersurface.
\end{conjecture}

 A published proof of this conjecture in dimensions less than or equal to seven can be found in~\cite{Eichmair:PET,EHLS,HuangLee}, building upon~\cite{SchoenYauPMT2,SchoenYauPMT1,SchoenYauPNAS}, and a preprint covering all dimensions is available as \cite{Lohkamp2} (with the borderline cases covered in \cite{HuangLee}).
Conjecturally, this result also follows in all dimensions by reduction from the time-symmetric case addressed in the preprint \cite{SchoenYau2017}.

In order to comply with the hypotheses of {Conjecture}~\ref{T12XII18.1}, from now on the manifold ${M}$ will be assumed to be the union of a compact set and of an asymptotic end diffeomorphic to $[0,\infty)\times {\mathbb S}^{n-1}$, where ${\mathbb S}^{k}$ denotes a $k$-dimensional sphere. In other words, ${M}$ will be assumed to
be a smooth, connected, non-compact  manifold with one end diffeomorphic to $[0,\infty)\times {\mathbb S}^{n-1}$.

\begin{Remark}
 \label{R9I19.1} {\rm
We note that the statement remains true with several asymptotically Euclidean ends, as large spheres in the remaining asymptotic ends provide barriers for the relevant constructions. (This shows in particular that orientability is not needed, as one can always pass to a double cover of the manifold; this introduces a second asymptotically Euclidean end which is innocuous for the problem at hand.)
Similarly outer trapped, or marginally outer-trapped compact boundaries are allowed; these will carry over to the hyperbolic setting as in Remark~\ref{R7I19.1} below.
}
\qed%
\end{Remark}
%
%We will show that {whenever 
%Conjecture~\ref{T12XII18.1} holds}, one can use the perturbation results in~\cite{CGNP} and the gluing constructions of~\cite{ChDelayExotic} to remove {some significant} restrictions in the proofs of positivity of asymptotically hyperbolic mass available elsewhere.

Let
$({M},g)$ denote an $n$-dimensional Riemannian manifold with the topology just-explained.
We use the ``exotic hyperbolic gluings'' of \cite{ChDelayExotic}, together with
the deformation results of~\cite{CGNP}, to establish our main result here, which is new for all dimensions unless the manifold is spin:

\begin{theorem}
  \label{T22VII18.3}
{Assume that Conjecture~\ref{T12XII18.1} is true in some dimension  $n\ge 3$}.
Consider an   asymptotically hyperbolic manifold $(M,g)$ 
 with spherical conformal infinity and  well-defined energy-momentum vector $\benmom$.
  If
%   \ptcr{here, and in other theorems, we do not say that this is conformally compact; and actually the new results by Yau and friends can be used, at least in low dimensions, to replace conformally compact by complete, should be said}
%
\bel{22VII18.1}
 R(g)\ge - n(n-1)
 \,,
\ee
then $\benmom$ is  causal  future directed or vanishes; in the last case $M$ is isometrically diffeomorphic to hyperbolic space.
\end{theorem}

\begin{Remark}
 \label{R7I19.1}
{\rm
In the asymptotically hyperbolic (AH) spin case, the timelike future-directed character of $\benmom$
has been established for manifolds which are complete with compact boundary, when the mean curvature $H$ of the boundary satisfies $H \le n-1$ (cf., e.g., \cite{ChHerzlich}).
Our method,
for reducing the Riemannian AH case to an asymptotically Euclidean (AE) one
(which  uses the replacement of $K$ by $K-g$, compare \eqref{10I19.1}-\eq{22VII18.4}),
transforms such boundaries in the AH regime to ones satisfying $H\le 0$ in the AE framework. This clarifies the somewhat mysterious AH condition $H \le n-1$.
\qed
}
\end{Remark}

\begin{Remark}
\label{r14IX19.1} {
\em
A precise set of asymptotic conditions needed for a well defined hyperbolic mass can be found e.g.\ in \cite[Proposition~2.2]{ChHerzlich}. It should be emphasised that these detailed asymptotic conditions are irrelevant for our proof in the following sense: given a metric with well defined energy-momentum, we start by deforming the metric in the asymptotic region so that the energy-momentum is almost unchanged but the perturbed metric is smoothly conformally compactifiable.  
When the Ricci scalar $R$ satisfies $R\le 0$ this can be done by repeating the arguments in \cite[Section~6]{ChDelayAH}. For metrics which fail to satisfy $R\le 0$ one can use the argument of \cite[Section~6]{ChDelayAH} together with a compactly supported deformation of the metric which compensates for the kernel of the linearised Yamabe operator, as in the proof of~\cite[Lemma~4.2]{DahlSakovich}. The deformed metric is then conformal to the hyperbolic metric in a neighborhood of the boundary, hence conformally smooth, with a smooth conformal factor by~\cite[Theorem~1.3]{ACF}. 
This allows one to carry-out subsequent arguments basing on smoothly compactifiable metrics with well defined mass aspect function, cf.\ \cite{CGNP}.
The fact that the glued metric used below to obtain a contradiction is likely to
%belong to a weighted Sobolev class near the boundary only, rather than a weighted H\"older one,
have less regularity at the conformal boundary does not affect the argument, as it can be slightly deformed again to a conformally smooth metric.
}
\qed
\end{Remark}

\nojdg{
\begin{Remark}
 \label{R18XII18.1a}
{\rm
Our theorems apply without further due to general relativistic data sets for which the energy of matter fields is positive and $\tr_g K$ vanishes, or in fact if
\begin{equation}\label{1I19.1}
\normp K\normp _g^2 - (\tr_g K)^2\ge 0
 \,.
\end{equation}
The same proofs establish the equivalents of Theorems~\ref{T22VII18.1} and  \ref{T22VII18.3} for general relativistic initial data $({M},g,K)$ satisfying the dominant energy condition in cases where the extrinsic curvature tensor $K$ is supported away   from  the asymptotically hyperbolic boundary at infinity regardless of \eqref{1I19.1}. A  version of Theorem~\ref{T22VII18.1} which allows non-compactly support curvature tensors $K$ would immediately lead to a corresponding generalisation of Theorem~\ref{T22VII18.3}.
}
\qed
\end{Remark}
}

The proof of  Theorem~\ref{T22VII18.3} proceeds by contradiction and uses perturbations of the metric. This approach prevents us from considering
the borderline cases. However, it is shown in \cite{LanHuangAHrigid} how to use our conclusion to
assert that $\benmom$ is timelike or vanishing, and vanishes only for hyperbolic space.

Let
$\mathbb H^{n }$ denote $n$-dimensional hyperbolic space (with sectional curvatures normalised to $-1$). A step in the proof of Theorem~\ref{T22VII18.3}, with some of interest of its own, is provided by the following:

\begin{thm}\label{rigid2}
 {Assume that Conjecture~\ref{T12XII18.1} is true in some dimension  $n\ge 3$}, 
and consider an $n$-dimensional   Riemannian manifold $({M},g)$ 
%\nojdg{either with
% \ptcr{ah and complete for the previous results; but conformally compactifiable here1}
%%
%\begin{enumerate}
% \item
% scalar curvature $R(g)$
%satisfying
%%
%$$
% R(g)\ge -n(n-1)
%  \,,
%$$
%%
%\item or equipped with a symmetric two-covariant tensor $K$ such that $(g,K)$ satisfies the dominant energy condition with cosmological constant $\Lambda = -n(n-1)/2$.
%%\erw{ok}
%\end{enumerate}
%}
%\onlyjdg{
with scalar curvature $R(g)$
satisfying
$$
 R(g)\ge -n(n-1)
  \,.
$$
If $({M},g)$ contains a region isometric to the complement of a compact set in $\mathbb H^{n }$,\nojdg{with $K$ vanishing there,}
then $({M},g)$ is  isometrically diffeomorphic to $\mathbb H^{n }$.
% \ptc{one could do a Maskit-type gluing at infinity along torii to get a sphere at infinity? but this would be wrong since torii can have negative mass?}
\end{thm}

We note that Theorem~\ref{rigid2}  
is already known to be true in dimensions $3\le n \le 7$ \cite{AnderssonGallowayCai}, where Conjecture~\ref{T12XII18.1} has already been established anyway, and
for spin manifolds~\cite{MinOo,AndDahl,Wang,ChHerzlich,LeeGR}).  
We do not assume that $M$ is spin, and our proof is  different from that of~\cite{AnderssonGallowayCai}.

We finish this introduction by noting that our analysis, together with Theorem~\ref{T12XII18.1} and an argument of Wang \cite{wang_uniqueness_2002}, implies uniqueness of higher-dimensional Anti de Sitter metric in the class of strictly static vacuum metrics with a negative cosmological constant and spherical conformal boundary at infinity.

\section{Definitions, notations and conventions}
 \label{sec:def}

For the convenience of the reader we summarise our definitions and conventions, which are {mostly} as in~\cite{ChDelayExotic}.

\subsection{Asymptotically hyperbolic manifolds}
 \label{ss20XI15.13}

We   model  hyperbolic space on the unit ball $B$  in $\R^n$, endowed with the metric
\begin{equation}\label{3XI18.4a}
 b=z^{-2}\delta
  \,,
\end{equation}
where $\delta$ is the Euclidean metric and $z(x)=\frac12(1-|x|_\delta^2)$.

Let ${(M,\zg)}$ be a smooth  $n$-dimensional
Riemannian manifold without  boundary such that, outside of a compact set, the manifold is diffeomorphic to $B$  minus a closed ball, and  there the metric $\zg$ equals $b$. The function $z$ is then extended to a smooth positive function on  $M$, still denoted by $z$.

Let $\varphi$  be a smooth strictly positive function
on $M$. For $k\in\N$ and $\alpha\in [0,1]$ we define
$C^{k,\alpha}_{\varphi}$ to be the space of $C^{k,\alpha}$
functions or tensor fields  for which the norm
% \erw{macro normp pour la norme ponctuelle, qui met une barre }
%
$$
\begin{array}{l}
\| u\| _{C^{k,\alpha}_{\varphi}(\ringg )}=\sup_{x\in
M}\sum_{i=0}^k\Big(
\normp \varphi  \mathring\nabla^{(i)}u(x)\normp _{\ringg }\\
 \hspace{3cm}+\sup_{0< d_{\ringg }(x,y)\le 1/2}\varphi(x)  \frac{\normp
\mathring\nabla^{(i)}u(x)-\mathring\nabla^{(i)}u(y)\normp _{\ringg }}{d^\alpha_{\ringg }(x,y)}\Big)
\end{array}
$$ is finite.
 Here  the norm $\normp\cdot \normp _{\ringg }$ and the covariant derivative $\znabla$ are defined using $\zg$, and $d _{\ringg }(x,y)$ is the $ \ringg  $-distance between $x$ and $y$.

A metric $g$ on $M$ will be called \emph{asymptotically hyperbolic}, or AH,  (with spherical conformal infinity) if  $g$ tends to $\zg$  in an appropriate function  space.
More precisely, given  a function space  $W$,
an asymptotically hyperbolic metric $g$ will be said to be of $M_{\zg+W} $ class if  $g-\zg\in W$.

In our context one is typically interested in $M_{\zg+W}$ metrics with $W=C^k_{z^{-\sigma}}$ with $\sigma>0$. In local coordinates as above such metrics decay to the model metric as $z^{\sigma-2}$, or as $z^{\sigma}$ in $\zg$-norm,  with derivatives satisfying uniform weighted estimates near the boundary. Further, there exists then a constant $C$ such that
\bel{gestimb}
  \normp g-\zg\normp _{\zg}+\normp \znabla g\normp _{\zg}+...+\normp \znabla^{(k)}g\normp _{\zg}\leq
 C z ^\sigma\,.
\ee
%$$
 For $\R\ni\sigma>k\ge 1$, $k\in \N$ and $g\in M_{\zg+C^k_{1,z^{-\sigma}}}$ the conformally  rescaled metric    $z^2 g$  can be extended  to the conformal boundary  $\mathbb {\mathbb S}^{n-1}$ of $M$, with the extension belonging to the  $C^{\lfloor\sigma\rfloor}$-differentiability class.

\subsection{Energy-momentum}

The definition of energy-momentum of an asymptotically hyperbolic   manifold requires choosing a coordinate system in which the metric manifestly approaches the hyperbolic metric. Let us denote by $\phi$ the choice of such a coordinate system. In the case of spherical conformal infinity,  which is of interest here,  such structures can be parameterised by the conformal group of ${\mathbb S}^{n-1}$. Indeed,
 let us denote by $\phi $ a coordinate system $(r,\theta^A)$, with $\theta^A$ local coordinates on ${\mathbb S}^{n-1}$, in which $g$ takes the form
\begin{equation}\label{5XI18.1}
  g
   \to_{r\to\infty}
 \frac { dr^2 }{1+r^2} + r^2 \underbrace{\ringh_{AB}(\theta) d\theta^A d\theta^B}_{=:\ringh} =: b
 \,,
\end{equation}
with $\ringh$ being the unit round metric on ${\mathbb S}^{n-1}$,
where the asymptotics is understood by requiring the $b$-norm $\normp g-b\normp _b $ of  $ g-b $ to decay to zero as one recedes to infinity. (To have a well defined energy-momentum one also needs fall-off rates, as well as  derivative-decay conditions~\cite{ChHerzlich,ChNagy}, which  will be implicitly assumed whenever relevant).  If $\hat \phi $ is another such coordinate system in which
\begin{equation}\label{5XI18.2}
  g
   \to_{\hat r\to\infty}
 \frac { d\hat r^2 }{1+\hat r^2} + \hat r^2 \underbrace{\ringh_{AB} (\hat \theta) d\hat \theta^A d\hat \theta^B}_{=:\hat{h }}
 \,,
\end{equation}
then there exists a conformal transformation $\Lambda $ of $({\mathbb S}^{n-1},\ringh)$ so that  $\hat{h}$ is obtained from  $\ringh$ by applying $\Lambda$ (cf., e.g., \cite{ChNagy,ChHerzlich}). We then write $\hat \phi  = \Lambda  \phi $,
hoping that the reader will not get confused by our use, later, of the same notation for the action of the conformal group on vectors by Lorentz transformations.

\nojdg{
We will write $(M ,
g ,K ,\phi)$ for an asymptotically hyperbolic initial data set $(M ,
g ,K )$ with an asymptotic structure $\phi$.
}

\onlyjdg{
We will write $(M , g ,\phi)$ for an asymptotically hyperbolic metric $(M ,
g )$ with an asymptotic structure $\phi$.
}

Recall that \emph{static Killing Initial Data} (KIDs) are solutions $V$ of the set of equations
\bea
 \zD_i \zD_j V
  =
   \Big(
   \mathrm{Ric}(b)_{ij}
     + n b_{ij}
 \Big)
        V
 \,,
  \label{29IV18.3}
\eea
where $b$ is the hyperbolic metric and $\zD$ its covariant derivative operator. This is an $(n+1)$-dimensional vector space equipped with a natural Lorentzian scalar product.

Letting $\{V_\mu\}_{\mu=0}^n$ be an orthonormal
basis of the space of static KIDs,
the energy-momentum vector $\benmom=(m_\mu)$ can be defined as
  \cite{HerzlichRicciMass} (compare \cite[Equation~(IV.40)]{BCHKK})
\begin{eqnarray}
 m_\mu
    &= &
   - % \frac{1}{16(n-2)\pi}
   \lim_{r\rightarrow\infty}\int_{{\mathbb S}^{n-1}(r) }  \red{D^j V_\mu}
    ( R{}^i{}_j - \frac {R{}}{n}\delta^i_j)
    \,
    d\sigma_i
     \,,
           \label{4XI18.6asdf}
    \end{eqnarray}
where $R_{ij}$ is the Ricci tensor of the metric $g$, $R$ its trace, and we have ignored an overall dimension-dependent positive multiplicative factor.

Using  the Poincar\'e-ball model, the ${\mathbb S}^{n-1}(r)$'s can be taken to be coordinate spheres of Euclidean radius $1-1/r$ accumulating at the conformal boundary ${\mathbb S}^{n-1}(+\infty)$  as $r$ tends to infinity.

The energy-momentum $\benmom$ is well defined   e.g.\ when $k\ge 1$, $\sigma>\frac n 2$, and $R+n(n-1)$, where $R$ is the curvature scalar of $g$, is in $L^1$; see~\cite{ChHerzlich}. In fact, rather weaker weighted-Sobolev  conditions suffice.

On the other hand, the mass aspect function is well defined if,
 e.g., the conformally rescaled metric $\tilde g:= z^{-n} g$  is $C^n(\overline M)$ in its own Gauss coordinates; this corresponds more or less to $\sigma\geq n$ and $k\ge n$.

\section{Localised ``Maskit'' gluings}
 \label{s4XI18.1}

The aim of this section is to analyse what happens with the energy-momentum under the localised  Maskit-type gluings of asymptotically hyperbolic manifolds of \cite[Section~3.5]{ChDelayExotic} for spherical conformal boundaries.  It is  clear that the analysis below can be adapted  to the  conformal gluings of
Isenberg, Lee and Stavrov~\cite{ILS} in the spherical case, and we will not discuss this case any further.
On the other hand it is not obvious how to adapt our calculations to non-spherical conformal boundaries, it would be of interest to settle this.
% \ptc{open problem}\erw{y reflechir un peu quand meme}

\nojdg{
The  gluing construction relevant for our purposes proceeds as follows:
Consider points $p_1$, $p_2$, lying on the conformal boundary of two asymptotically hyperbolic initial data sets
$(M_1,g_1,K_1)$ and $(M_2,g_2,K_2)$  satisfying the dominant energy condition, with total energy-momentum vectors
$ \benmom^{\myone} \equiv (\enmom_{\mu}^{\myone})$ and $\benmom ^{\mytwo } \equiv (\enmom_{\mu}^{\mytwo})$. We will assume that both $(M_1, g_2)$ and $(M_2, g_2)$ have spherical conformal infinity (by this we mean that the conformal class of the boundary metric is that of the canonical metric on the sphere),
with  the extrinsic curvature tensors asymptoting to zero, and with well-defined total energy-momentum.
As shown in \cite{ChDelayExotic}
for deformations of data sets preserving the vacuum condition or for scalar curvature deformations preserving an inequality, or in Appendix~\ref{s10XII18.1} below for deformations of data sets preserving the dominant energy condition,
 for all $\varepsilon>0$ sufficiently small we can construct new initial data sets
$(M_i, g_{i,\red{\varepsilon} },K_{i,\red{\varepsilon} })$, $i=1,2$,
such that

\begin{enumerate}
\item
the metrics coincide with the hyperbolic metric in coordinate half-balls $\mcU_\varepsilon^{\myone }$ of radius $\varepsilon$ (measured with respect the ``compactified'' metric) around $p_1$ and
$\mcU_\varepsilon^{\mytwo }$ around $p_2$, and the $K_{i,\red{\varepsilon} }$'s vanish there;
\item
the new initial data sets coincide with the original ones outside of balls of radius $2\varepsilon$;
\item the total energy-momenta of the new metrics tend to the original ones as $\varepsilon$ tends to zero.
\end{enumerate}
}

\onlyjdg{
The  gluing construction relevant for our purposes proceeds as follows:
Consider points $p_1$, $p_2$, lying each on the conformal boundary of two asymptotically hyperbolic manifolds
$(M_1,g_1)$ and $(M_2,g_2)$, each metric with scalar curvature larger than or equal $-n(n-1)$, with total energy-momentum vectors
$ \benmom^{\myone} \equiv (\enmom_{\mu}^{\myone})$ and $\benmom ^{\mytwo } \equiv (\enmom_{\mu}^{\mytwo})$. We will assume that both $(M_1, g_2)$ and $(M_2, g_2)$ have spherical conformal infinity (by this we mean that the conformal class of the boundary metric is that of the canonical metric on the sphere),
and have a  well-defined total energy-momentum.
As shown in \cite{ChDelayExotic}
 for all $\varepsilon>0$ sufficiently small we can construct  new metrics
$(M_i, g_{i,\red{\varepsilon} })$, $i=1,2$,  each with scalar curvature larger than or equal $-n(n-1)$,
such that:

\begin{enumerate}
\item
the metrics coincide with the hyperbolic metric in coordinate half-balls $\mcU_\varepsilon^{\myone }$ of radius $\varepsilon$ (as measured with respect to the ``compactified'' metric) around $p_1$ and
$\mcU_\varepsilon^{\mytwo }$ around $p_2$;
\item  the new metrics coincide with the original ones outside of balls of radius $2\varepsilon$;

\item the total energy-momenta of the new metrics tend to the original ones as $\varepsilon$ tends to zero.
\end{enumerate}
}

Choose some asymptotic structures $\phi_i$ on $(M_i,g_i)$ and let
$$
 \benmom ^{i,\varepsilon} \equiv (\enmom_{\mu}^{i,\varepsilon})
$$
denote the energy-momenta of
\onlyjdg{
$(M_i, g_{i,\red{\varepsilon} },\phi_i)$.
}
\nojdg{
$(M_i, g_{i,\red{\varepsilon} },K_{i,\red{\varepsilon} },\phi_i)$.
}
 As already mentioned,  the construction in~\cite{ChDelayExotic} guarantees that
\begin{equation}\label{4XI18.3}
  \enmom_{\mu}^{i,\varepsilon} \to_{\varepsilon\to 0} \enmom_{\mu}^{i }
  \,,
  \ i=1,2
  \,.
\end{equation}
\begin{remark}
 \label{R10XI18.1}
{\rm
For our further purposes we need to complement \eqref{4XI18.3} with a decay rate, which can be read-off from~\cite[Section 3.3]{ChDelayExotic}.
We note that  the background hyperbolic metric, denoted in the current work by $b$,
 is  denoted by $\ringg $ there
(recall that $\ringg  $  is used here to denote a metric on $M$ which coincides with $b$ on the asymptotic end but not necessarily everywhere, e.g.\ because of a different topology of $M$),
and that we use here  $s$  instead of the symbol $b$ used  in~\cite{ChDelayExotic} for the
rate of radial decay of the glued metrics.
Let us thus assume that
$$
 \mbox{$\normp g-b\normp_b=O(r^{-\sigma}) $ with  $\sigma> (n-1)/2+s$}
$$
for some  $s\in[n/2 , (n+1)/2)$.
By \cite[Remark~3.6]{ChDelayExotic}
 and the comments at the end of   Section 3.3 in~\cite{ChDelayExotic}, the $b$-norm of the correction $h_{\varepsilon}$ to the metric arising from the gluing procedure is then of order $o (\varepsilon^{\sigma-s}r^{-s})$.
Now recall that the mass is the limit, as $r$ tends to infinity, of an integral on a sphere of radius $r$ of a quantity, say $\mathbb U$ (cf. \eq{4XI18.6asdf}).
In this limit only the part of $\mathbb U$ linear in
$e_{\varepsilon}=g_{\varepsilon}-b=g-b+h_{\varepsilon}
=e+h_{\varepsilon}$ %=e+O(\varepsilon^{\sigma-s}r^{-s})$
matters.
 This implies  (cf.\ the last, unnumbered, equation in the proof of \cite[Theorem~3.7]{ChDelayExotic})
\begin{equation}\label{estimasse}
  \enmom_{\mu}^{i,\varepsilon}-\enmom_{\mu}^{i }=o(\varepsilon^{\sigma-s})
  \,,
  \ i=1,2
  \,.
\end{equation}
When working in weighted H\"older spaces, the decay rate  $s>n/2$ is needed for the final mass to be well defined, while one needs $\sigma\leq n$ for the initial mass to be non-zero. The resulting best  approach rate ${\sigma-s}$   is  as close as desired to, but smaller than, $n/2$.
In our positivity proof below, the decay rate $\sigma=n$ is obtained by invoking the density results of~\cite{DahlSakovich}, independently of the  decay rate of the metric under consideration. We will see that this suffices for the positivity theorem in space dimensions $n\ge 5$.

A more careful analysis, which provides
\begin{equation}
  \enmom_{\mu}^{i,\varepsilon}-\enmom_{\mu}^{i }=
 o(\varepsilon^{n/2})
  \,,
  \ i=1,2
  \,,
   \label{14XII18.1}
\end{equation}
is needed when $n=4$. To obtain this estimate one notes that all the gluing results in~\cite{ChDelayExotic} can be traced back to Theorem~3.3 there, \nojdg{or an initial data version thereof,} where the solution of the problem at hand is constructed in weighted Sobolev spaces. In view of the already mentioned results in~\cite{DahlSakovich}, in \cite[Theorem~3.3]{ChDelayExotic} we can take $\sigma=n$, and choose $b$ there (which coincides with $s$ here) to be $n/2$
(the mass still being well defined, the correction being in $H^{k+2}_{1,z^{-n/2}}$ in the notation of
\cite{ChDelayExotic}). This provides the improved decay rate \eqref{14XII18.1}; compare~\cite[Remark~3.4]{ChDelayExotic}.
}
\qed
\end{remark}

Let
$$
 \mcH \subset  \Hyp^n\subset \R^n
$$
denote the
equatorial hyperplane
in the Poincar\'e ball model. So $\mcH$ is a totally geodesic hypersurface which ``cuts the hyperbolic space in half'', and its conformal completion
intersects the conformal boundary of $\Hyp^n$ at the equator. We denote by $\Hyp^n_+$ the part of $\Hyp^n $ which lies ``above'' $\mcH$ and by $\Hyp^n_-$ the part of $\Hyp^n $ which lies ``below'' $\mcH$.
More precisely, $\Hyp^n_+$ is the part of $\Hyp^n$ with the height coordinate larger than or equal zero, and $\Hyp^n_-$ is the part of $\Hyp^n$ with the height coordinate smaller than or equal zero.
(While this is completely irrelevant at this stage, we note that for typographical reasons, related to the calculations in Equation~\eqref{9XI18.2} and following, it will be convenient to think of the height as being determined by
 the  first coordinate
in $\R^n$, thus the equator on the conformal boundary will be defined as the intersection of the unit coordinate sphere in $\R^n$ with the hyperplane $\{x^1=0\}$.)

Recall that the connected component of the group of conformal isometries of a sphere ${\mathbb S}^{n-1}$ is isomorphic to the component of the identity of the Lorentz group in dimension $n+1$ and the action,
which we denote by   $\Lambda \benmom$,
 of
an element $\Lambda$ of the conformal group of the sphere
on  the energy-momentum vector $\benmom$, is the standard linear action of the Lorentz group on $\R^{1,n}$.

\begin{figure}
    \centering
\begin{tabular}{m{3.3cm} m{.3cm} m{3.3cm} m{.3cm} m{2.5cm} }
     \includegraphics[width=.25\textwidth,angle=-90]{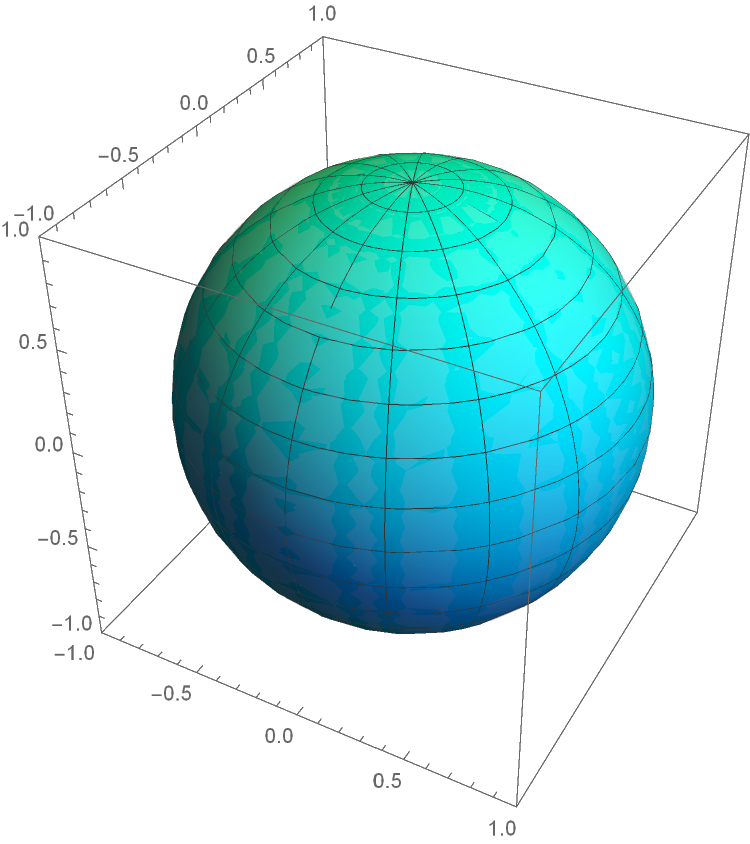}
    &
     $\overset{\mathrm{glue}}{\longrightarrow}$
    &
   \includegraphics[width=.25\textwidth,angle=-90]{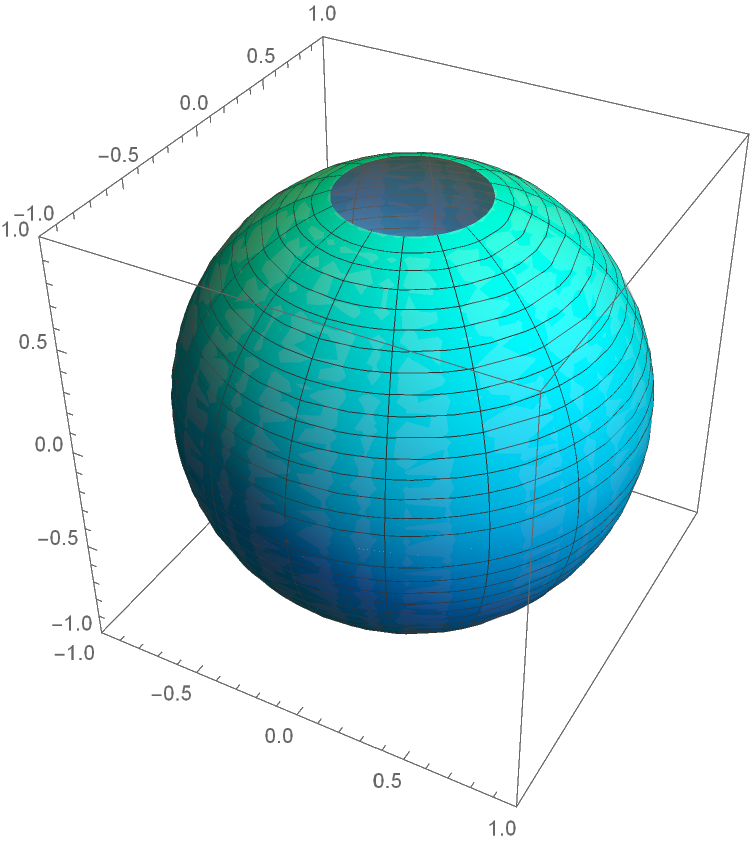}
    &
     $\overset{\Lambda_\varepsilon^{\myone }}{\longrightarrow}$
    &
    \includegraphics[width=.25\textwidth,angle=-90]{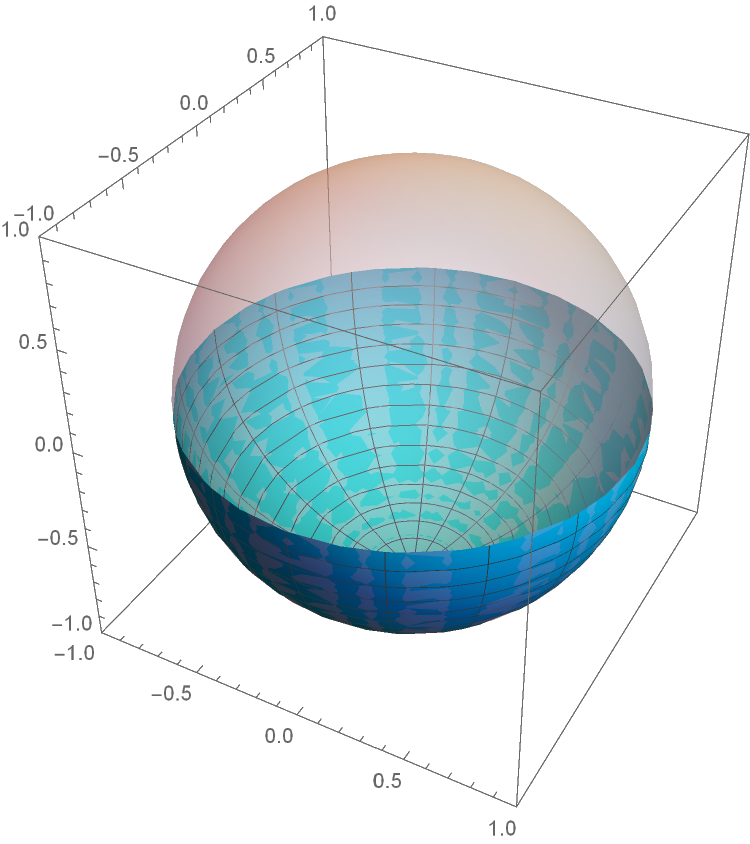}
    \end{tabular}
 \caption{\label{F8XI18.1} A small spherical cap $D^1_\varepsilon\subset {\mathbb S}^{n-1}$   (middle picture) is mapped to the upper-half sphere by the conformal transformation $\Lambda_\varepsilon^{\myone }$. The final metric coincides with the hyperbolic one \red{in a region which includes} the upper half-ball. We assume the topology of conformal infinity to be spherical, but the interior does not have to be topologically trivial. Thus, the topology inside the final upper half-sphere is that of a half-ball, with the boundary of the lower-half sphere (including the equatorial hyperplane) bounding the remaining topology of $M_1$.}
\end{figure}

Let $D_\varepsilon^{\myi }$ denote the $(n-1)$-dimensional ball centered at $p_i$, contained in the conformal boundary, obtained by intersecting $\mcU_\varepsilon^{\myi }$ with the conformal boundary; see the middle figure in Figure~\ref{F8XI18.1}.
Let us denote by $\Lambda_\varepsilon^{\myone }$ the conformal transformation of the conformal boundary which maps $D_\varepsilon^{\myone }$ to the \emph{upper half-sphere}.
See the right figure in~Figure~\ref{F8XI18.1}.
%
%
%\begin{figure}
%    \centering
%    \begin{tabular}{m{3.3cm} m{.5cm} m{3.3cm} m{.5cm} m{2.5cm} }
%   \includegraphics[width=.25\textwidth]{WholeSphere}
%    &
%     $\overset{\mathrm{glue}}{\longrightarrow}$
%    &
%   \includegraphics[width=.25\textwidth]{jarmulka2}
%    &
%     $\overset{\Lambda_\varepsilon^{\myone }}{\longrightarrow}$
%    &
%    \includegraphics[width=.25\textwidth]{jarmulka4}
%    \end{tabular}
% \caption{\label{F8XI18.1} A small spherical cap $D^1_\varepsilon\subset {\mathbb S}^{n-1}$   (middle picture) is mapped to the upper-half sphere by the conformal transformation $\Lambda_\varepsilon^{\myone }$. The final metric coincides with the hyperbolic one throughout the upper half-ball. We assume the topology of conformal infinity to be spherical, but the interior does not have to be topologically trivial. Thus, the topology inside the final upper half-sphere is that of a half-ball, with the boundary of the lower-half sphere (including the equatorial hyperplane) bounding the remaining topology of $M_1$.}
%\end{figure}
%

In physics terminology, $\Lambda_\varepsilon^{\myone }$ is a boost in the direction of $-p_1$,
 when ${\mathbb S}^{n-1}$ is viewed as a coordinate sphere in $\R^n$; this can be seen e.g. from \cite[Section~3.4]{CGNP}. The velocity parameter of the boost tends to the speed of light as $\varepsilon$ tends to zero.

After applying $\Lambda_\varepsilon^{\myone }$ we obtain a  set
\nojdg{$(M_1, g_{1,\red{\varepsilon} },K_{1,\red{\varepsilon} },  \Lambda_\varepsilon^{\myone }\phi_1)$}
\onlyjdg{$(M_1, g_{1,\red{\varepsilon} },  \Lambda_\varepsilon^{\myone }\phi_1)$}  with energy-momentum vector $\Lambda_\varepsilon^{\myone }\benmom ^{1,\varepsilon}$.   $M_1$ contains a region,   denoted by $\mcV_\varepsilon^{\myone } $,
defined as
the ``half-hyperbolic-space'' $\Hyp^n_+$  in the coordinate system $\Lambda_\varepsilon^{\myone }\phi_1$.
 The boundary $\partial \mcV_\varepsilon^{\myone } $ corresponds in the local coordinates to $\mcH$, is totally geodesic, and the metric $ g_{1,\red{\varepsilon} }$ is exactly hyperbolic throughout  $\Hyp^n_+$ and in a neighborhood of $\mcH \approx \partial\mcV_\varepsilon^{\myone }  $.
 Thus
\begin{equation}\label{7XI18.1}
 \Hyp^n_+ \approx  \mcV_\varepsilon^{\myone } \subset \mcU_\varepsilon^{\myone }  \subset M_1
  \,,
\end{equation}
where $\approx$ in \eq{7XI18.1} means ``isometrically diffeomorphic to''.

Letting, as before, $\{V_\mu\}_{\mu=0}^n$ be an orthonormal
basis of the space of static KIDs, we have
\begin{eqnarray}
 (\Lambda_\varepsilon^{\myone }\benmom ^{1,\varepsilon})_\mu
    &= &
   - % \frac{1}{16(n-2)\pi}
   \lim_{r\rightarrow\infty}\int_{{\mathbb S}^{n-1}(r) }  \red{D^j V_\mu}
    ( R{}^i{}_j - \frac {R{}}{n}\delta^i_j)
    \,
    d\sigma_i
     \,,
           \label{4XI18.6}
\end{eqnarray}
where now $R_{ij}$ is the Ricci tensor of the metric $g_{1,\red{\varepsilon} } $.

Let ${\mathbb S}^{n-1}_-(r)$ denote the part of the sphere which lies under the equator, and ${\mathbb S}^{n-1}_+(r)$ the part above the equator. Since the trace-free part of the Ricci tensor of $g_{1,\red{\varepsilon} } $ vanishes on $\Hyp^n_+$,  only the integrals on ${\mathbb S}^{n-1}_-(r)$ contribute:
\begin{eqnarray}
 (\Lambda_\varepsilon^{\myone }\benmom ^{1,\varepsilon})_\mu
    &= &
   - % \frac{1}{16(n-2)\pi}
   \lim_{r\rightarrow\infty}\int_{{\mathbb S}^{n-1}_-(r) }  \red{D^j V_\mu}
    ( R{}^i{}_j - \frac {R{}}{n}\delta^i_j)
     \,
    d\sigma_i
     \,.
           \label{4XI18.7 }
\end{eqnarray}

\begin{figure}
    \centering
\begin{tabular}{m{3.3cm} m{.3cm} m{3.3cm} m{.3cm} m{2.5cm} }
   \includegraphics[width=.25\textwidth,angle=-90]{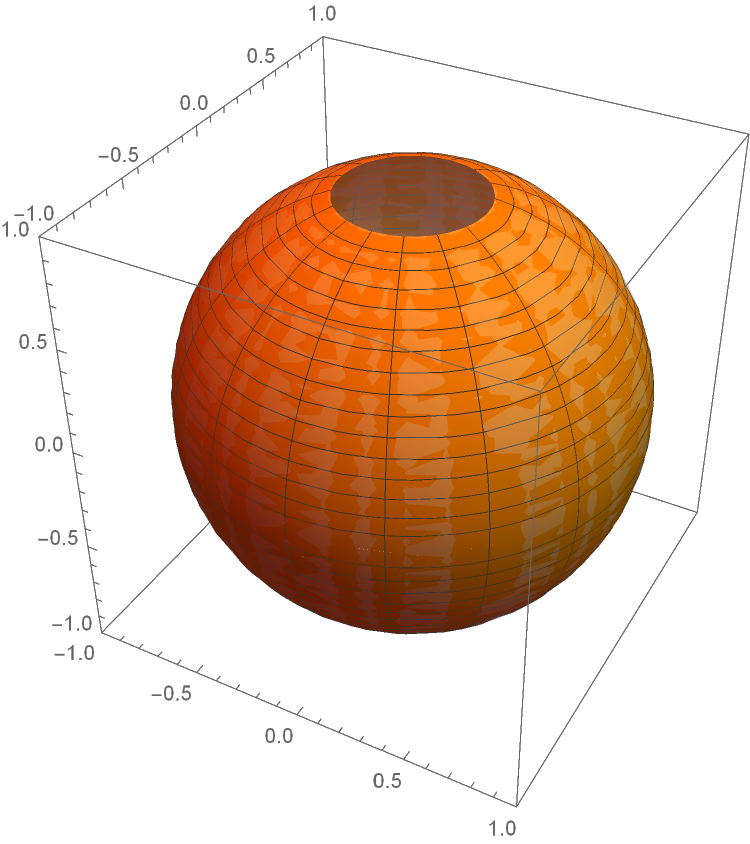}
    &
     $\overset{\Lambda_\varepsilon^{\mytwo }}{\longrightarrow}$
    &
   \includegraphics[width=.25\textwidth,angle=-90]{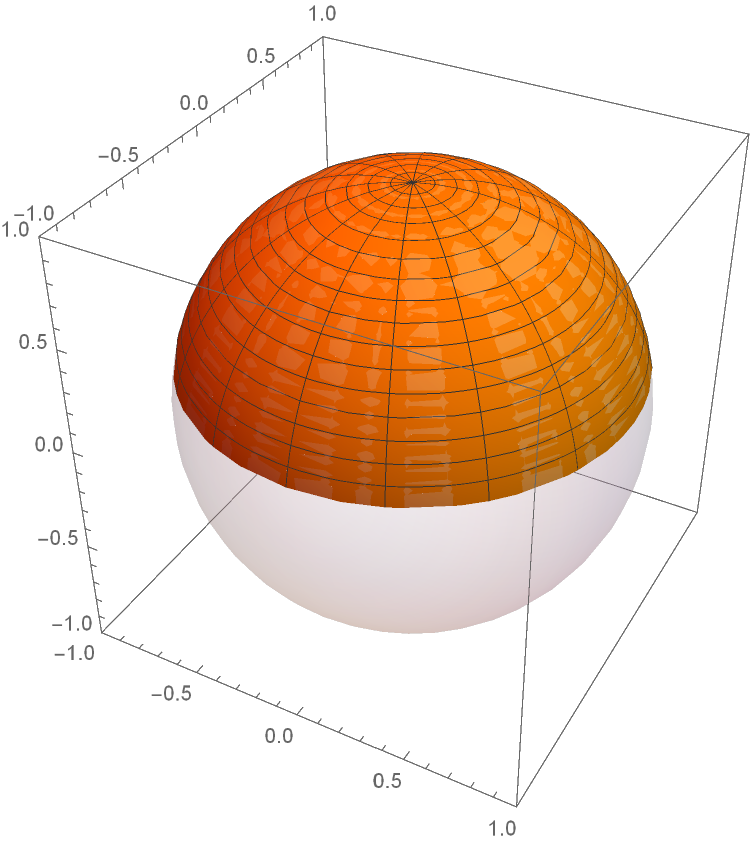}
    &
     $\overset{\mathrm{glue}}{\longrightarrow}$
    &
    \includegraphics[width=.25\textwidth,angle=-90]{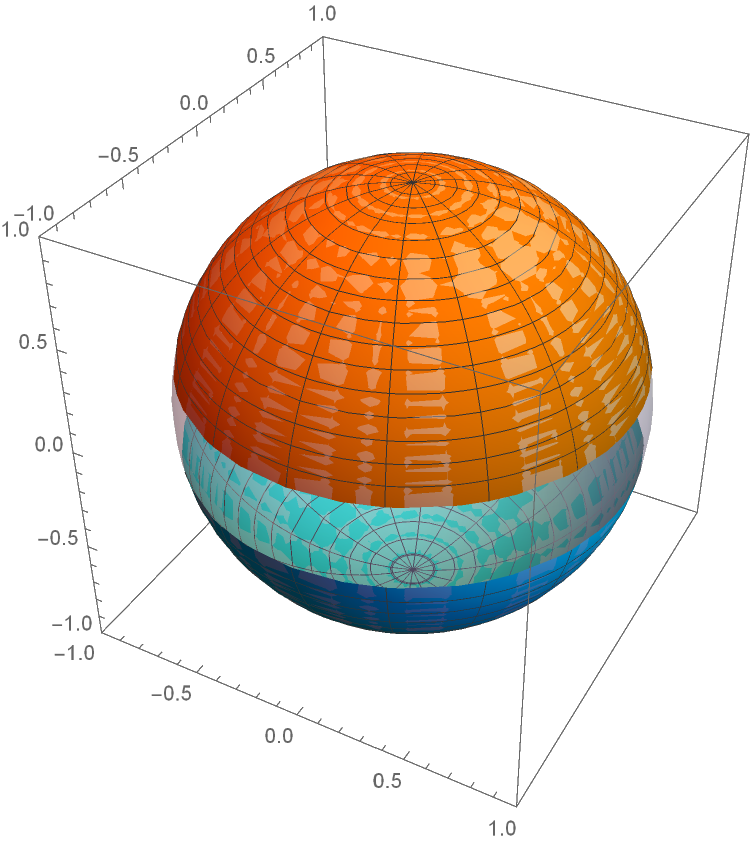}
    \end{tabular}
 \caption{\label{F12XI18.1} A small spherical cap $D^2_\varepsilon\subset {\mathbb S}^{n-1}$   (left picture) is mapped to the lower-half sphere by the conformal transformation $\Lambda_\varepsilon^{\mytwo }$. The final manifold $(M,g_\varepsilon)$ is obtained by gluing together the manifold from the  right Figure~\ref{F8XI18.1} with the manifold from the middle figure here. The metric  $g_\varepsilon$ coincides with $g_{1,\varepsilon}$ on the ``lower half'' of the final manifold and with $g_{2,\varepsilon}$ on the ``upper half'', except for two transition strips near the equatorial hyperplane, \red{and is exactly hyperbolic in an overlapping strip enclosing the equatorial hyperplane shown in the right figure}.}
\end{figure}

Similarly there exists a conformal transformation of the conformal boundary, which will be denoted by  $\Lambda_\varepsilon^{\mytwo }$, which maps $D_\varepsilon^{\mytwo }$ to the \emph{lower half-sphere}.
See Figure~\ref{F12XI18.1}.
After applying this transformation, $M_2$ contains a region,
which we denote by $ \mcV_\varepsilon^{\mytwo }$,
isometrically diffeomorphic to the  half-hyperbolic-space $\mcH_-$, as well as its   totally geodesic boundary $\mcH$, with the metric being exactly hyperbolic in a neighborhood of $\mcH$, and throughout $ \mcV_\varepsilon^{\mytwo }$. Analogously to \eqref{7XI18.1} we have
$$
 \Hyp^n_- \approx  \mcV_\varepsilon^{\mytwo }  \subset\mcU_\varepsilon^{\mytwo }  \subset M_2
  \,.
$$
One also finds
\begin{eqnarray}
 (\Lambda_\varepsilon^{\mytwo }\benmom ^{2,\varepsilon})_\mu
    &= &
   - % \frac{1}{16(n-2)\pi}
   \lim_{r\rightarrow\infty}\int_{{\mathbb S}^{n-1} (r) }  \red{D^j V_\mu}
    ( R{}^i{}_j - \frac {R{}}{n}\delta^i_j)
     \,
    d\sigma_i
     \nonumber
    \\
    &= &
   - % \frac{1}{16(n-2)\pi}
   \lim_{r\rightarrow\infty}\int_{{\mathbb S}^{n-1}_+(r) }  \red{D^j V_\mu}
    ( R{}^i{}_j - \frac {R{}}{n}\delta^i_j)
     \,
    d\sigma_i
     \,,
                \label{4XI18.11}
\end{eqnarray}
where $R_{ij}$ is now the Ricci tensor of the metric $g_{2,\red{\varepsilon} } $.

Since the metric equals the hyperbolic metric near $\mcH$ in both manifolds,\nojdg{ and the extrinsic curvature tensor vanishes  near $\mcH$ in both manifolds, we can construct a new family of asymptotically hyperbolic initial data sets $(M ,g_{ \red{\varepsilon} },K_{ \red{\varepsilon} })$ by gluing together $M_1\setminus  \mcV_\varepsilon^{\myone }$ with $M_2\setminus  \mcV_\varepsilon^{\mytwo }$ along $\mcH$, with the obvious asymptotic structure $\phi$  and the obvious smooth initial data set  $( g_{ \red{\varepsilon} },K_{ \red{\varepsilon} })$ .} we can construct a new family of AH metrics $(M ,g_{ \red{\varepsilon} } )$ by gluing together $M_1\setminus  \mcV_\varepsilon^{\myone }$ with $M_2\setminus  \mcV_\varepsilon^{\mytwo }$ along $\mcH$, with the obvious asymptotic structure $\phi$  and the obvious smooth metric  $g_{ \red{\varepsilon} }$ .

Denote by $ \benmom ^{ \varepsilon} \equiv (  \enmom ^{ \varepsilon} _\mu)$ the total energy-momentum vector of\nojdg{ $(M ,g_{ \red{\varepsilon} },K_{ \red{\varepsilon} },\phi)$.}\onlyjdg{ $(M ,g_{ \red{\varepsilon} },\phi)$.} We have
\begin{eqnarray}
 \enmom ^{ \varepsilon} _\mu
    &= &
   - % \frac{1}{16(n-2)\pi}
   \lim_{r\rightarrow\infty}\int_{{\mathbb S}^{n-1} (r) }  \red{D^j V_\mu}
    ( R{}^i{}_j - \frac {R{}}{n}\delta^i_j)
    \,
    d\sigma_i
    \nonumber
    \\
    &= &
   - % \frac{1}{16(n-2)\pi}
   \lim_{r\rightarrow\infty}\int_{{\mathbb S}^{n-1}_-(r) }  \red{D^j V_\mu}
    ( R{}^i{}_j - \frac {R{}}{n}\delta^i_j)
    \,
    d\sigma_i
    \nonumber
    \\
    &&
   - % \frac{1}{16(n-2)\pi}
   \lim_{r\rightarrow\infty}\int_{{\mathbb S}^{n-1}_+(r) }  \red{D^j V_\mu}
    ( R{}^i{}_j - \frac {R{}}{n}\delta^i_j)
    \,
    d\sigma_i
    \nonumber
    \\
     & = &
 (\Lambda_\varepsilon^{\myone }\enmom ^{1,\varepsilon})_\mu
  +
 (\Lambda_\varepsilon^{\mytwo }\enmom ^{2,\varepsilon})_\mu
     \,.
           \label{4XI18.6a}
\end{eqnarray}
In other words,
\begin{eqnarray}
 \benmom ^{ \varepsilon}
    &= &
 \Lambda_\varepsilon^{\myone }\benmom ^{1,\varepsilon}
  +
 \Lambda_\varepsilon^{\mytwo }\benmom ^{2,\varepsilon}
     \,.
           \label{4XI18.6b}
\end{eqnarray}

Thus, the energy-momentum is additive under the above Maskit-type gluing,
in a sense made precise by \eqref{4XI18.6b}.

\section{Positivity for asymptotically hyperbolic Riemannian manifolds}
 \label{s29X18.1}

 For the convenience of the reader we recall  Conjecture~\ref{T12XII18.1}:

\begin{conjecture}
\label{T30X18.1}
Let  $(M,g,K )$ be an $n$-dimensional general relativistic initial data set, $n\ge 3$, satisfying the dominant energy condition such that
$g$ is flat outside of a compact set and $K$ vanishes there.
 Then $(M,g)$ can be isometrically embedded in Minkowski space-time so that $K$ is the second fundamental form of the embedded hypersurface.
\end{conjecture}
 
\newcommand{\shouldbenothing}{}
\newcommand{\shouldbeone}{{\red{1}}}
\newcommand{\shouldbetwo}{{\red{2}}}
\newcommand{\shouldbethree}{{\red{3}}}

We first deduce the following 

\begin{thm}\label{rigid1}
Consider an $n$-dimensional   Riemannian manifold $({M},g)$, $n \ge 3$,  
with scalar curvature $R(g)$
satisfying
$$
 R(g)\ge -n(n-1)
  \,.
$$
Suppose that Conjecture \ref{T30X18.1} holds in this dimension.
If $({M},g)$ contains a region isometric to the complement of a compact set in $\mathbb H^{n }$,\nojdg{with $K$ vanishing there,}
then $({M},g)$ is  isometrically diffeomorphic to $\mathbb H^{n }$.
\end{thm}

\noindent{\sc Proof}
Consider an asymptotically hyperbolic manifold  {$(M,g\shouldbenothing)$} as in the statement of the theorem. 
The pair $(M,g\shouldbenothing)$ can be viewed as a general relativistic initial data set $(M,g\shouldbenothing,K\shouldbenothing\equiv 0)$ for Einstein equations with negative cosmological constant $\Lambda = - n(n-1)/2$, with matter fields with positive energy density $\rho\shouldbenothing:= R(g\shouldbenothing) + n(n-1) $, and with vanishing matter current $J\shouldbenothing:= 0$.
 But one can also view $(M,g\shouldbenothing)$ as an initial data set
\begin{equation}\label{10I19.1}
 (M,g_\shouldbeone:=g\shouldbenothing,K_\shouldbeone:= {-}g_\shouldbeone)
\end{equation}
for Einstein equations with vanishing cosmological constant, with positive energy density $\rho_\shouldbeone=\rho\shouldbenothing$, and with vanishing matter current $J_\shouldbeone\equiv 0$. This follows immediately from the following trivial calculation
\begin{equation}\label{22VII18.4}
  R(g_\shouldbeone) \equiv  R(g\shouldbenothing) =   \rho + 2 \Lambda =   \rho + \normp K_\shouldbeone\normp ^2_{g_\shouldbeone} - \tr_{g_\shouldbeone} (K_\shouldbeone)^2
   \,.
\end{equation}
Let us denote by $\mcU\subset M$ the region of $M$ where $g_\shouldbeone\equiv g\shouldbenothing$ is exactly hyperbolic. Then $(\mcU,g_\shouldbeone)$ can be embedded into $(n+1)$-dimensional Minkowski space-time $\R^{1,n}$ as a hyperboloid so that its second fundamental form coincides with $K_\shouldbeone$. By passing to a subset of $\mcU$ if necessary there exists $R>0$ so that $\mcU$ can be identified with
\begin{equation}\label{22VII18.5}
  \mcU = \{t={-}\sqrt {1 + |\vec x|^2}\,, \
   |\vec x| > R
    \}
    \,.
\end{equation}
Let $f$ be any smooth function which coincides with ${-}\sqrt {1 + |\vec x|^2}$ for $R\le |\vec x| \le R+1$, which is constant for large $|\vec x|$, and such that the graph of $f$ over  the set $\{|\vec x| \ge R\}$ is spacelike. Let $M_\shouldbetwo$ denote this graph, let $g_\shouldbetwo$ be the metric induced by the Minkowski metric on the graph of $f$, and let $K_\shouldbetwo$ be the second fundamental form of the graph. Let $(M_\shouldbethree,g_\shouldbethree,K_\shouldbethree)$ be obtained by taking the union of
$$(M\setminus \mcU,g\shouldbenothing,-g\shouldbenothing)\equiv(M\setminus \mcU,g_\shouldbeone,-g_\shouldbeone)$$
with $(M_\shouldbetwo,g_\shouldbetwo,K_\shouldbetwo)$ and identifying points at the joint boundary $|\vec x |=R$.
Then $(M_\shouldbethree,g_\shouldbethree,K_\shouldbethree)$ is a smooth general relativistic initial data satisfying the dominant energy condition which is exactly flat at large distances.
By {Conjecture}~\ref{T30X18.1} the data set $(M_\shouldbethree,g_\shouldbethree,K_\shouldbethree)$ arises from Minkowski space-time, in particular $M_\shouldbethree$ is spin. %In particular $(M_\shouldbethree,g_\shouldbethree)$ is spin, and t
The result follows now from\onlyjdg{ \cite{ChHerzlich,Wang}.}\nojdg{\cite{ChHerzlich,Wang} in the Riemannian case, or \cite{CJL} in the initial data setting.}
%\xout{So $(M,g)$ is a spin manifold with timelike past-pointing total energy-momentum, which is not possible by~\cite{Wang,ChHerzlich}.}
%
\qedskip

We also note the following result:% 
\footnote{The hypotheses have been chosen for simplicity, the result holds under much more general asymptotic conditions, compare Remark~\ref{r14IX19.1}.} 

\begin{theorem} \label{T22VII18.2a}
{Assume that Conjecture~\ref{T30X18.1} is true in a dimension  $n\ge 3$}.
Let $({M},g)$,  be an $n$-dimensional
manifold with scalar curvature $R[g] \ge -n(n-1)$ with a metric which is smoothly conformally compactifiable and with spherical conformal infinity.
Suppose that 
%there exists $\epsilon>0$ such that 
in the coordinates of \eq{5XI18.1} the metric $g$ approaches the hyperbolic metric $\mathring g$ as $O(r^{-n})$, 
%$O(r^{-n/2-\epsilon})$, 
%with  $|\mathring D g|_{\mathring g} =  O(r^{-n/2-\epsilon})$,
and assume that $R[g]+n(n-1)=O(r^{-n-1})$.
Then the total energy-momentum vector of $({M},g)$ \emph{cannot} be timelike past-pointing.
\end{theorem}

\proof
Since every three dimensional manifold is spin,  this is the usual AH Positive Energy Theorem\onlyjdg{~\cite{Wang,ChHerzlich}}\nojdg{~\cite{CJL}} when $n=3$.

Suppose, thus, that $n\ge 4$ and that the claim is wrong. 
By hypothesis the metric is conformally smooth, so that the
% decay with  $|\mathring D g|_{\mathring g} =  O(r^{-n/2-\epsilon})$,
%One can make a small deformation of the metric, which does not change the causal character of the energy-momentum vector, 
%so that the 
asymptotic conditions of \cite{CGNP} are met.
Using a deformation 
of the metric as in \cite[Corollary~1.4]{CGNP}  (one can alternatively invoke \cite[Theorem~1.2]{ChDelayAH} if $n\ge 8$), one obtains an AH metric $g_1$   with
\bel{22VII18.2-}
 R(g_1)\ge - n(n-1) \  \mbox{ and } \  R(g_1)\not \equiv - n(n-1)
 \,,
\ee
such that $g_1$ has negative mass aspect function.
As explained in~\cite{AnderssonGallowayCai},   one can now construct on $M$ a new metric $g_2 $ which satisfies
\bel{22VII18.2}
 R(g_2)\ge - n(n-1) \  \mbox{ and } \  R(g_2)\not \equiv - n(n-1)
 \,,
\ee
and which coincides with the hyperbolic metric $\hypg $
in a neighborhood of the conformal boundary at infinity. This contradicts 
Theorem~\ref{rigid1}. 
\qed

We are ready now to pass to the

\medskip

\noindent {\sc Proof of Theorem~\ref{T22VII18.3}:}
\nojdg{For further reference we will carry out the argument as far as possible  (and thus not as far as the final conclusions) for general relativistic initial data sets satisfying the dominant energy condition, to make it clear how any generalisation of Theorem~\ref{T22VII18.1} to include the extrinsic curvature tensor $K$ would provide a generalisation of Theorem~\ref{T22VII18.3}.

Suppose that there exists an asymptotically hyperbolic initial data set
$(M_1,g_1,K_1)$
}\onlyjdg{
Suppose that there exists an AH manifold
$(M_1,g_1)$}
with spacelike or null past-pointing energy-momentum vector $\benmom^1$ with respect to an asymptotic structure $\phi_1$.
\nojdg{If $K_1\equiv 0$, w}We can make a  deformation of the metric as in \cite[Proposition~6.2]{ChDelayAH}  to achieve that $\normp g-\ringg\normp _\ringg$
  is $O(r^{-n})$ in the coordinates of \eq{5XI18.1} while maintaining the space-like or past-directed causal character of the energy-momentum vector.\nojdg{
In the general case, by approximation results in the spirit of those in~\cite{DahlSakovich} one should similarly be able to obtain that $\normp g-b\normp _b$   and $\normp K\normp _b$
are $O(r^{-n})$ in the coordinates of \eq{5XI18.1}; alternatively one could just assume these fall-off rates.}

Changing $\phi_1$ (i.e., applying a conformal transformation at the conformal boundary) if necessary, we can choose the asymptotic structure $\phi_1$ so that
\begin{equation}\label{8XI18.1}
  \benmom^{1} \equiv (\enmom^{1}_0, \vec \enmom^{1}) \in \R\times \R^n
\end{equation}
satisfies
\begin{equation}\label{8XI18.2}
   \enmom^{1}_0 < 0
   \,.
\end{equation}
 This follows from the fact that the connected component of the identity of the Lorentz group acts transitively on the set of spacelike vectors with given length.

Choose $p_1\in \mathbb S^{n-1}$ to be $-\vec \enmom^1$.

\onlyjdg{Let $
 (M_1, g_{1,\red{\varepsilon} } ,\phi_1)$
be a  family  as in Section~\ref{s4XI18.1}, constructed by a gluing procedure which, in half-plane-model coordinates, replaces $g_1$ by the hyperbolic metric in a half-ball of radius $\varepsilon$ centred at $p_1$.
Set
$$
 (M_2, g_{2,\red{\varepsilon} } ,\phi_2)
 :=
 (M_1, g_{1,\red{\varepsilon} } ,\phi_1)
  \,.
$$
 }
\nojdg{Let $
 (M_1, g_{1,\red{\varepsilon} },K_{1,\red{\varepsilon} },\phi_1)$
be a  family  as in Section~\ref{s4XI18.1}, constructed by a gluing procedure which, in half-plane-model coordinates, replaces $g_1$ by the hyperbolic metric and $K_1$ by zero in a half-ball of radius $\varepsilon$ centred at $p_1$.
Set
$$
 (M_2, g_{2,\red{\varepsilon} },K_{2,\red{\varepsilon} },\phi_2)
 :=
 (M_1, g_{1,\red{\varepsilon} },K_{1,\red{\varepsilon} },\phi_1)
  \,.
$$
}

 Let $(M,g_\varepsilon,\phi)$, with energy-momentum $\benmom^\varepsilon$, be constructed as in Section~\ref{s4XI18.1} by gluing together
$ (M_1, g_{1,\red{\varepsilon} } ,\phi_1)$ and
$ (M_2, g_{2,\red{\varepsilon} } ,\phi_2)$.

The maps $\Lambda_\varepsilon^{\myone }$ of Section~\ref{s4XI18.1} act on $\benmom^1_\varepsilon$ as boosts in the direction of $\vec \enmom^1 $.
Let $R_\varepsilon$ be a rotation of the sphere by $\pi$ around
some arbitrarily chosen axis passing through the equator,  thus $R_\varepsilon$ maps the north pole to the south pole.
If we align the coordinate system so that $\vec \enmom^1$ points in the direction of the first coordinate axis, we can write
\begin{equation}\label{9XI18.2}
  \Lambda_\varepsilon^{\myone }
  =
   \left(\begin{array}{cccc}
      \gamma_\varepsilon &-\gamma_\varepsilon v_\varepsilon & 0 & 0   \\
      -\gamma_\varepsilon v_\varepsilon & \gamma_\varepsilon & 0 & 0  \\
      0 & 0 & 1 & 0   \\
      0 & 0 & 0    & \Id
    \end{array}
    \right)
  \,,
  \quad
   R_\varepsilon = \left(\begin{array}{cccc}
      1 & 0 & 0 & 0    \\
     0 & -1 & 0 & 0  \\
      0 & 0 & -1 &  0   \\
      0 & 0 & 0    & \Id
    \end{array}
    \right)
  \,,
\end{equation}
with $v_\varepsilon\in (0,1)$, $\gamma_\varepsilon:= (1-v_\varepsilon^2)^{-1/2}$, and $\Id$ being the $(n-2)\times (n-2)$ identity matrix.

Let us choose $\Lambda_\varepsilon^{\mytwo}$ to be $R_\varepsilon \Lambda_\varepsilon^{\myone }$.
Then the space-part of $\Lambda^2_{\varepsilon} \benmom^{2}=\Lambda^2_{\varepsilon} \benmom^{1}$ cancels out that of $\Lambda^1_{\varepsilon} \benmom^{1}$:
\begin{equation}\label{4XII18.13}
  \Lambda^1_{\varepsilon} \benmom^{1} + \Lambda^2_{\varepsilon} \benmom^{1}=
  \Lambda^1_{\varepsilon} \benmom^{1} +R_\varepsilon \Lambda^1_{\varepsilon} \benmom^{1}=
  %\Lambda^1_{\varepsilon}
  2 \gamma_\varepsilon (\underbrace{\enmom^{1}_0}_{<0} - v_\varepsilon \underbrace{\enmom^1_1}_{> |\enmom^{1}_0| >0}, \vec 0)
  \,.
\end{equation}
Hence,
\begin{eqnarray}
 \benmom ^{ \varepsilon}
    &= &
 \Lambda_\varepsilon^{\myone }  \benmom ^{1,\varepsilon}
  +
 \Lambda_\varepsilon^{\mytwo }\benmom ^{2,\varepsilon}
 =
 \Lambda_\varepsilon^{\myone }  \benmom ^{1,\varepsilon}
  +
 R_\varepsilon\Lambda_\varepsilon^{\myone }\benmom ^{1,\varepsilon}
 \nonumber
\\
    &= &
    (
 \Lambda_\varepsilon^{\myone }
 +
 R_\varepsilon  \Lambda_\varepsilon^{\myone }  )(
 \benmom ^{1}+
 \benmom ^{1,\varepsilon}
 - \benmom ^{1}
 )
 \nonumber
\\
    &= &
   \gamma_\varepsilon\bigg(
    2(\underbrace{\enmom^{1}_0-v_\varepsilon \enmom^1_1}_{<0}, \vec 0)
 +
 \nonumber
\\
 &&
      \gamma_\varepsilon^{-1} \Lambda^1_{\varepsilon}
    \Big(
  \underbrace{
  \big(1 + (\Lambda_\varepsilon^{\myone } )^{-1} R_\varepsilon
 \Lambda_\varepsilon^{\myone }
 \big)
 ( \benmom ^{1,\varepsilon} - \benmom^1)
 }_{=: (*)}
 \Big)
  \bigg)
     \,.
           \label{4XI18.14}
\end{eqnarray}
Keeping in mind that $ \gamma_\varepsilon^{-1} \Lambda^1_{\varepsilon}$ has bounded entries, we wish to show that $(*)$ goes to zero as $\varepsilon$ goes to zero.
For this, we need to analyse the  matrix $(\Lambda_\varepsilon^{\myone } )^{-1} R_\varepsilon
 \Lambda_\varepsilon^{\myone }$. Since $(\Lambda_\varepsilon^{\myone } )^{-1}$ coincides with $\Lambda_\varepsilon^{\myone } $ with $v_\varepsilon$ changed to its negative,
  we find
\begin{equation}\label{9XI18.3}
  (\Lambda_\varepsilon^{\myone } )^{-1} R_\varepsilon
 \Lambda_\varepsilon^{\myone }
  =
   \left(\begin{array}{cccc}
      \gamma_\varepsilon^2(1+v_\varepsilon^2) &-2\gamma_\varepsilon^2 v_\varepsilon & 0 & 0   \\
      2\gamma_\varepsilon^2 v_\varepsilon & -\gamma_\varepsilon^2(1+v_\varepsilon^2) & 0 & 0  \\
      0 & 0 & -1 & 0   \\
      0 & 0 & 0    & \Id
    \end{array}
    \right)
  \,.
\end{equation}
If we view the sphere ${\mathbb S}^{n-1}$ as embedded in  $\R^n$ with Euclidean coordinates $(x^i, x^n)$, with $-\vec \enmom^1$ the north pole and $\vec \enmom^1$ the south pole, then the Lorentz transformation $\Lambda^\myone_\varepsilon$ of \eqref{9XI18.2}
%\erw{ajout detail  correspondance  avec la ref CGNP entre parenthese et proposition changement equa dessous \\ -- \\ ptc: je vais reverifier, mais d'autre part, si on repete les notations de \cite[Section~3.4]{CGNP} ici, autant les utiliser partout; donc soit on enleve, soit on commence a changer partout....}
corresponds to the conformal transformation (cf., e.g., \cite[Section~3.4]{CGNP})
% with the corresponding notations $v_\varepsilon=-\frac{\sinh\alpha}{\cosh\alpha}$, $e$ the  unit vector of $\R^n$ with first component $1$, so $M_{\alpha,e}$ correspond to $\Lambda_\varepsilon $ and is related to the conformal transformation $\Phi_{\alpha,e}^\infty$ there,  denoted $\Lambda_\epsilon$ here. We recall that when $\Lambda_\epsilon$ act on a vector $(y_0,y)$ of $\R\times\R^n$, then the corresponding conformal map on $\mathbb S^{n-1}$ transform his space direction $x=y/|y|$ to the space direction of its image)
%%
%\begin{equation}\label{10XI18.1}
%{\mathbb S}^{n-1} \ni (x^i, x^n) \mapsto \frac 1 {\gamma_\varepsilon(1-v_\varepsilon x^n)}( x^i, \gamma_\varepsilon(x^n-v_\varepsilon)) \in {\mathbb S}^{n-1}
% \,.
%\end{equation}
%%
%\erw{remplacer par :}
%
\begin{equation}\label{10XI18.1b}
{\mathbb S}^{n-1} \ni (x^1,x^i) \mapsto \frac 1 {\gamma_\varepsilon(1-v_\varepsilon x^1)}(\gamma_\varepsilon(x^1-v_\varepsilon),x^i) \in {\mathbb S}^{n-1}
 \,.
\end{equation}
If we denote by $ \theta$ the angle off the axis passing through the south pole and the north pole, then $v_\varepsilon$ has to be chosen so that $\varepsilon$ is mapped to $\pi/2$. Equivalently, if $n^i$ is a unit vector in Euclidean $\R^{n-1}$, Equation~\eqref{10XI18.1b} should give
%%
%\begin{equation}\label{10XI18.3}
%{\mathbb S}^{n-1} \ni (\sin(\varepsilon) n^i, \cos(\varepsilon) ) \mapsto \frac 1 {\gamma_\varepsilon(1-v_\varepsilon \cos(\varepsilon))}( n^i, 0) \in {\mathbb S}^{n-1}
% \,,
%\end{equation}
%%
%\erw{remplacer par :}
%%
\begin{equation}\label{10XI18.3b}
{\mathbb S}^{n-1} \ni ( \cos(\varepsilon),\sin(\varepsilon) n^i ) \mapsto \frac 1 {\gamma_\varepsilon(1-v_\varepsilon \cos(\varepsilon))}(0,\sin(\varepsilon) n^i) \in {\mathbb S}^{n-1}
 \,,
\end{equation}
i.e., $v_\varepsilon= \cos(\varepsilon)$, hence
\begin{equation}\label{10XI18.2}
  \gamma_\varepsilon = \frac 1 {\sin(\varepsilon)} =\varepsilon^{-1}
   \left( 1 + O(\varepsilon)\right)
\end{equation}
for small $\varepsilon$. It follows that the term $ (*)$ in \eqref{4XI18.14} can be estimated as
\begin{eqnarray}
 |(*)|
 &
  \equiv &
  |\big(1 + (\Lambda_\varepsilon^{\myone } )^{-1} R_\varepsilon
 \Lambda_\varepsilon^{\myone }
 \big)
 ( \benmom^{1,\varepsilon} - \benmom^1)
 |
 \nonumber
\\
 &
 \le &
 C(1 +  \varepsilon^{-2}) | \benmom ^{1,\varepsilon} - \benmom^1 | \le C^2 \varepsilon^{ \frac n 2 -2}
     \,,
           \label{4XI18.14a}
\end{eqnarray}
where we used \eqref{estimasse}. So $ (*)$ will tend to zero if $n \ge 5$.
The case $n=4$ follows in the same way using \eqref{14XII18.1}.
This implies that for $\varepsilon$ small enough the energy-momentum vector $\benmom^\varepsilon$ of $(M,g_\varepsilon,\phi)$
is
timelike past pointing. 
One can perturb the metric slightly 
if needed (compare Remark~\ref{r14IX19.1}) so that the asymptotic conditions  of Theorem~\ref{T22VII18.2a} are satisfied, with the energy-momentum remaining timelike past pointing, 
obtaining a 
contradiction with that theorem. 
We conclude that $\benmom$ is causal future-pointing, or vanishes.
%
%As already pointed out, one can now use~\cite{LanHuangAHrigid}  to show that $\benmom$ is  timelike or vanishing, and vanishes only for hyperbolic space.
%
\qed

\bigskip

\noindent{\sc Acknowledgements:} We thank Michael Eichmair for useful discussions related to Theorem \ref{T12XII18.1},
and Greg Galloway for many useful comments.
PTC was supported in part by
the Austrian Science Fund (FWF) under project  P29517-N27  and by
the Polish National Center of Science (NCN) 2016/21/B/ST1/00940.
ED was supported by the grant ANR-23-CE40-0010  (project EINSTEIN-PPF)
and the grant ANR-24-CE40-0702 (project OrbiScar)  of the French National Research Agency ANR.

\providecommand{\bysame}{\leavevmode\hbox to3em{\hrulefill}\thinspace}
\providecommand{\MR}{\relax\ifhmode\unskip\space\fi MR }
\providecommand{\MRhref}[2]{%
  \href{http://www.ams.org/mathscinet-getitem?mr=#1}{#2}
}
\providecommand{\href}[2]{#2}


\begin{thebibliography}{10}

\bibitem{AbbottDeser}
L.F. Abbott and S.~Deser, \emph{Stability of gravity with a cosmological
  constant}, Nucl.\ Phys. \textbf{B195} (1982), 76--96.

\bibitem{AnderssonGallowayCai}
L.~Andersson, M.~Cai, and G.J. Galloway, \emph{Rigidity and positivity of mass
  for asymptotically hyperbolic manifolds}, Ann.\ H.~Poincar\'e \textbf{9}
  (2008), 1--33, arXiv:math.dg/0703259. \MR{MR2389888 (2009e:53054)}

\bibitem{ACF}
L.~Andersson, P.T. Chru\'{s}ciel, and H.~Friedrich, \emph{On the regularity of
  solutions to the {Y}amabe equation and the existence of smooth hyperboloidal
  initial data for {E}insteins field equations}, Commun.\ Math.\ Phys.
  \textbf{149} (1992), 587--612. \MR{MR1186044 (93i:53040)}

\bibitem{AndDahl}
L.~Andersson and M.~Dahl, \emph{Scalar curvature rigidity for asymptotically
  locally hyperbolic manifolds}, Annals of Global Anal.\ and Geom. \textbf{16}
  (1998), 1--27, arXiv:dg-ga/9707017.

\bibitem{BCHKK}
H.~Barzegar, P.T. Chru\'{s}ciel, and M.~{H\"orzinger}, \emph{Energy in
  higher-dimensional spacetimes}, Phys.\ Rev.\ D \textbf{96} (2017), 124002, 25
  pp., arXiv:1708.03122 [gr-qc].

\bibitem{ChDelayAH}
P.T. Chru\'{s}ciel and E.~Delay, \emph{Gluing constructions for asymptotically
  hyperbolic manifolds with constant scalar curvature}, Commun.\ Anal.\ Geom.
  \textbf{17} (2009), 343--381, arXiv:0711.1557[gr-qc]. \MR{2520913
  (2011a:53052)}

\bibitem{ChDelayExotic}
\bysame, \emph{Exotic hyperbolic gluings}, Jour.\ Diff.\ Geom. \textbf{108}
  (2018), 243--293, arXiv:1511.07858 [gr-qc].

\bibitem{CGNP}
P.T. Chru\'{s}ciel, G.J. Galloway, L.~Nguyen, and T.-T. Paetz, \emph{{On the
  mass-aspect function and positive energy theorems for asymptotically
  hyperbolic manifolds}}, Class. Quantum Grav. \textbf{35} (2018), 115015,
  arXiv:1801.03442 [gr-qc].

\bibitem{ChHerzlich}
P.T. Chru\'{s}ciel and M.~Herzlich, \emph{The mass of asymptotically hyperbolic
  {R}iemannian manifolds}, Pacific Jour.\ Math. \textbf{212} (2003), 231--264,
  arXiv:math/0110035 [math.DG]. \MR{MR2038048 (2005d:53052)}

\bibitem{CJL}
P.T. Chru\'{s}ciel, J.~Jezierski, and S.~\L\c{e}ski, \emph{The{ Trautman-Bondi}
  mass of hyperboloidal initial data sets}, Adv.\ Theor.\ Math.\ Phys.
  \textbf{8} (2004), 83--139, arXiv:gr-qc/0307109. \MR{MR2086675 (2005j:83027)}

\bibitem{ChNagy}
P.T. Chru\'{s}ciel and G.~Nagy, \emph{The {H}amiltonian mass of asymptotically
  {anti -- de Sitter} spacetimes}, Class.\ Quantum Grav. \textbf{18} (2001),
  L61--L68, arXiv:hep-th/0011270.

\bibitem{ChNagyATMP}
\bysame, \emph{The mass of spacelike hypersurfaces in asymptotically {anti --
  de Sitter} spacetimes}, Adv.\ Theor.\ Math.\ Phys. \textbf{5} (2001),
  697--754, arXiv:gr-qc/0110014.

\bibitem{ChruscielSimon}
P.T. Chru\'{s}ciel and W.~Simon, \emph{Towards the classification of static
  vacuum spacetimes with negative cosmological constant}, Jour.\ Math.\ Phys.
  \textbf{42} (2001), 1779--1817, arXiv:gr-qc/0004032.

\bibitem{DahlSakovich}
M.~Dahl and A.~Sakovich, \emph{{A density theorem for asymptotically hyperbolic
  initial data satisfying the dominant energy condition}}, Pure Appl.\ Math.\
  Q. \textbf{17} (2015), 1669--1710, arXiv:1502.07487 [math.DG]. \MR{4376092}

\bibitem{Eichmair:PET}
M.~Eichmair, \emph{{The Jang equation reduction of the spacetime positive
  energy theorem in dimensions less than eight}}, Commun.\ Math.\ Phys.
  \textbf{319} (2013), 575--593, arXiv:1206.2553 [math.dg]. \MR{3040369}

\bibitem{EHLS}
M.~Eichmair, L.-H. Huang, D.A. Lee, and R.~Schoen, \emph{The spacetime positive
  mass theorem in dimensions less than eight}, Jour.\ Eur.\ Math.\ Soc.\ (JEMS)
  \textbf{18} (2016), 83--121, arXiv:1110.2087 [math.DG]. \MR{3438380}

\bibitem{GHHP}
G.W. Gibbons, S.W. Hawking, G.T. Horowitz, and M.J. Perry, \emph{Positive mass
  theorem for black holes}, Commun.\ Math.\ Phys. \textbf{88} (1983), 295--308.

\bibitem{HerzlichRicciMass}
M.~Herzlich, \emph{Computing asymptotic invariants with the {R}icci tensor on
  asymptotically flat and asymptotically hyperbolic manifolds}, Ann.\ Henri
  Poincar\'e \textbf{17} (2016), 3605--3617, arXiv:1503.00508 [math.DG].
  \MR{3568027}

\bibitem{LanHuangAHrigid}
L.-H. Huang, H.C. Jang, and D.~Martin, \emph{Mass rigidity for hyperbolic
  manifolds}, Commun.\ Math.\ Phys. \textbf{376} (2019), 2329--2349,
  arXiv:1904.12010 [math.DG].

\bibitem{HuangLee}
L.-H. Huang and D.A. Lee, \emph{Equality in the spacetime positive mass
  theorem}, Commun.\ Math.\ Phys. \textbf{376} (2020), 2379--2407,
  arXiv:1706.03732 [math.DG]. \MR{4104553}

\bibitem{ILS}
J.~Isenberg, J.M. Lee, and I.~Stavrov~Allen, \emph{Asymptotic gluing of
  asymptotically hyperbolic solutions to the {E}instein constraint equations},
  Ann.\ Henri Poincar\'e \textbf{11} (2010), 881--927, arXiv:0910.1875
  [math-dg]. \MR{2736526 (2012h:58037)}

\bibitem{LeeGR}
D.A. Lee, \emph{Geometric relativity}, Graduate Studies in Mathematics, vol.
  201, AMS, Providence, RI, 2017.

\bibitem{Lohkamp2}
J.~{Lohkamp}, \emph{{The Higher Dimensional Positive Mass Theorem II}},
  (2016), arXiv:1612.07505 [math.DG].

\bibitem{Maerten}
D.~Maerten, \emph{Positive energy-momentum theorem in asymptotically {anti-de
  Sitter} spacetimes}, Ann.\ H.Poincar\'e \textbf{7} (2006), 975--1011,
  arXiv:math.DG/0506061. \MR{MR2254757 (2007d:83016)}

\bibitem{MinOo}
M.~Min-Oo, \emph{Scalar curvature rigidity of asymptotically hyperbolic spin
  manifolds}, Math.\ Ann. \textbf{285} (1989), 527--539.

\bibitem{SchoenYauPNAS}
R.~Schoen and S.-T. Yau, \emph{Complete manifolds with nonnegative scalar
  curvature and the positive action conjecture in general relativity}, Proc.\
  Nat.\ Acad\. Sci.\ U.S.A. \textbf{76} (1979), 1024--1025. \MR{MR524327
  (80k:58034)}

\bibitem{SchoenYauPMT1}
\bysame, \emph{On the proof of the positive mass conjecture in general
  relativity}, Commun.\ Math.\ Phys. \textbf{65} (1979), 45--76. \MR{MR526976
  (80j:83024)}

\bibitem{SchoenYauPMT2}
\bysame, \emph{Proof of the positive mass theorem. {II}}, Commun.\ Math.\ Phys.
  \textbf{79} (1981), 231--260. \MR{612249}

\bibitem{SchoenYau2017}
\bysame, \emph{Positive scalar curvature and minimal hypersurface
  singularities},   Surveys in Diff. Geom. \textbf{24} (1) (2017), arXiv:1704.05490 [math.DG].

\bibitem{Wang}
X.~Wang, \emph{Mass for asymptotically hyperbolic manifolds}, Jour.\ Diff.\
  Geom. \textbf{57} (2001), 273--299. \MR{MR1879228 (2003c:53044)}

\bibitem{wang_uniqueness_2002}
X.~Wang, \emph{On the uniqueness of the {A}d{S} spacetime}, Acta Math.\ Sin.
  (Engl. Ser.) \textbf{21} (2005), 917--922 (en), arXiv:math/0210165.
  \MR{2156971}

\end{thebibliography}
\end{document}